\documentclass[a4paper, 10pt]{article}
\usepackage{amsmath,amsthm,amsfonts,amssymb,color}
\usepackage{a4wide}
\usepackage[dvips]{graphicx}

\theoremstyle{plain}
\newtheorem{thm}{Theorem}[section]

\newtheorem{lem}[thm]{Lemma}
\newtheorem{teigi}[thm]{Definition}

\newtheorem{rem}[thm]{Remark}

\newcommand{\ds}{\displaystyle}

\newcommand{\bC}{\mathbb{C}}
\newcommand{\bE}{\mathbb{E}}

\newcommand{\cF}{\mathcal{F}}
\newcommand{\cI}{\mathcal{I}}

\newcommand{\cK}{\mathcal{K}}

\newcommand{\cO}{\mathcal{O}}
\newcommand{\cP}{\mathcal{P}}

\newcommand{\bS}{\mathbb{S}}

\newcommand{\R}{\mathbb{R}}

\newcommand{\cX}{\mathcal{X}}
\newcommand{\x}{\boldsymbol{x}}
\newcommand{\y}{\boldsymbol{y}}

\newcommand{\1}{\boldsymbol{1}}

\newcommand{\diam}{\text{diam}}

\author{Masatake Hirao\footnote{
Department of Information and Science Technology, 
Aichi Prefectural University, 
1522-3 Ibaragabasama, Nagakute, Aichi 480-1198 Japan
(e-mail: hirao@ist.aichi-pu.ac.jp)
}}

\title{Finite frames, frame potentials and determinantal point processes on the sphere}
\begin{document}

\maketitle




\begin{abstract}
Herein,  
we address the expectations of frame potentials of three types of \textit{determinantal point processes} (DPPs) on the unit sphere $\bS^d$:
(i) spherical ensembles on $\bS^2$; (ii) harmonic ensembles on $\bS^d$ and (iii) jittered sampling point processes on $\bS^d$. 
The random point configurations generated by such DPPs converge more rapidly towards 
\textit{finite unit norm tight frames} (FUNTFs) than the Poisson point processes on the unit sphere.
\end{abstract}





\section{Introduction}
\label{sect:intro}
Let $\bS^d = \{ \x = (x_1, \ldots, x_{d+1}) \in \R^{d+1} \mid \| \x \| = 1 \}$ be the $d$-dimensional unit sphere, 
where $\| \x \| = \sqrt{\langle \x, \x \rangle}$ for the Euclidean norm and $\langle \x, \y \rangle = x_1 y_1 + \cdots + x_{d+ 1} y_{d + 1}$ for the standard inner product of $\R^{d+1}$.
We are interested in the ^^ ^^ goodness" of random point configurations generated by 
\textit{determinantal point processes} (DPPs)
on the unit sphere $\bS^d$, which are used in, e.g., a fermion model in quantum mechanics.

Herein, we deal with the \textrm{frame potential} presented previously by Benedetto and Fickus~\cite{BF03}. 
Let $\{ \x_i \}_{i = 1}^{n}$  be an $n$-point subset of $\bS^d$. 
The frame potential of $\{ \x_i \}_{i = 1}^{n}$ is 
\begin{equation}
\label{eq:FP}
\textrm{FP}(\{ \x_i \}_{i = 1}^{n}) = \sum_{i = 1}^{n} \sum_{j = 1}^{n} |\langle \x_i, \x_j \rangle|^2.
\end{equation}

The frame potential concept was introduced as a measure of 
well-distribution 
of point configurations on the unit sphere.
Further, its minimizers were characterized, e.g., 
the vertices of a regular polygon or a regular polyhedron inscribed in $\bS^1$ or $\bS^2$, 
respectively, always attain the lower bound of (\ref{eq:FP}).

The frame potential is closely related to finite frames. 
If $n \geq d + 1$, the finite subset $\{ \x_i \}_{i = 1}^{n}$ of $\bS^d$ that attains the lower bound of~(\ref{eq:FP}) is a 
\textit{finite unit norm tight frame} (FUNTF) 
for $\R^{d+1}$. 
The subset $\{ \x_i \}_{i = 1}^{n}$ of $\R^{d+1}$ is called a \textit{finite frame for $\R^{d+1}$}
if there are two constants $0 < \alpha \leq \beta < \infty$ such that 
$\alpha \| \x \|^2 \leq  \sum_{i = 1}^{n} |\langle \x, \x_i \rangle|^2 \leq \beta \| \x \|^2$
for all $\x \in \R^{d+1}$.
$\{ \x_i \}_{i = 1}^{n}$ is called a \textit{finite tight frame for $\R^{d+1}$} if $\alpha = \beta$.
If a finite tight frame $\{ \x_i \}_{i = 1}^{n}$ is a subset of $\bS^{d}$, 
it can be called the FUNTF for $\R^{d+1}$.

This study aims to calculate the expectations of the frame potentials of spherical ensembles on $\bS^2$, 
harmonic ensembles on $\bS^d$ and jittered sampling point processes on $\bS^d$; see Theorem~\ref{thm:main1}. 
For comparison with these results, note that, 
if $\{ \x_i \}_{i = 1}^{n}$ denotes an $n$-point Poisson point process on $\bS^d$, 
we can easily check that 
$\bE (\mathrm{FP}(\{ \x_{i} \}_{i = 1}^n)) = n^2/(d + 1) + d n/(d + 1)$.
Further, we aim to show that the random point configurations generated by such DPPs converge more rapidly towards FUNTFs than the Poisson point processes on the unit sphere; 
see Theorem~\ref{thm:main2} and Remark~\ref{rem:1}~(i).
\begin{thm}
\label{thm:main1}
\begin{itemize}
\item[(i)]
If $\{ \x_i \}_{i = 1}^{n}$ denotes an $n$-point spherical ensemble on $\bS^2$, 
then 
\begin{equation*}
\bE (\mathrm{FP}(\{ \x_{i} \}_{i = 1}^n) = \frac{n^2}{3} + \frac{4 n^2}{(n + 1)(n+2)}.
\end{equation*}
\item[(ii)] 
If $\{ \x_i \}_{i = 1}^{n}$ denotes an $n$-point harmonic ensemble on $\bS^d$ 
with $n = \dim (\cP_L (\bS^d)) = \binom{d + L}{d} + \binom{d+L-1}{d}$, 
then 
\begin{align*}
\bE (\mathrm{FP}(\{ \x_{i} \}_{i = 1}^n))
&=  \frac{n^2}{d+1} + \left( \frac{d}{L +d} - \frac{d (d - 1)}{(d + L) (d + 2 L + 1) (d + 2 L -1)}\right) \binom{d + L}{d} 
\\
& \left ( =  \frac{n^2}{d+1} +  \cO(n^{\frac{d - 1}{d}} )\quad (n \rightarrow \infty) \right ).
\end{align*}
\item[(iii)] If $\{ \x_i \}_{i = 1}^{n}$ denotes the $n$-point jittered sampling point process on $\bS^d$
given by Definition~\ref{def:jspp}, then
\begin{align*}
\frac{n^2}{d + 1} \leq \bE(\mathrm{FP} (\{ \x_i \}^n_{i = 1}))
\leq \frac{n^2}{d + 1}  + n - n \left (1 - \frac{c^2}{2 n^{2/d}} \right )^2.
\end{align*}
\end{itemize}
\end{thm}

The frame properties can be often expressed based on the \textit{analysis matrix} $X_n$ defined as 
$X_n \x = (\langle \x, \x_1 \rangle, \ldots, \langle \x, \x_n \rangle)^{T} \in \R^n$ for all $\x \in \R^{d + 1}$
and its transpose \textit{synthesis matrix} $X_n^{T}$.
Let $\cI_{d + 1}$ be the identity matrix of dimension $d + 1$ and 
$\| \cdot \|_{\cF}$ be the Frobenius norm of the matrix.

If $\{ \x_i \}_{i = 1}^{n}$ is a FUNTF for $\R^{d+1}$, 
then it satisfies $n^{-1} X_{n} X_{n}^T = (d + 1)^{-1} \cI_{d + 1}$ and vice versa (cf., Lemmas 2.2 and 2.3 of Ehler~\cite{E11}). 
Theorem~\ref{thm:main2} indicates that the DPPs on the unit sphere converges towards FUNTFs in terms of the expectation of the squared Frobenius norm. 
\begin{thm}
\label{thm:main2}
\begin{itemize}
\item[(i)] If $\{ \x_i \}_{i = 1}^{n}$ denotes an $n$-point spherical ensemble on $\bS^2$, 
then 
\[
\bE (\| \frac{1}{n} X_{n}^{T} X_{n} - \frac{1}{3} \cI_{3} \|_{\cF}^{2} ) = \frac{4}{(n + 1)(n + 2)}.
\]
\item[(ii)] If $\{ \x_i \}_{i = 1}^{n}$ denotes an $n$-point harmonic ensemble on $\bS^d$ with $n = \dim (\cP_L (\bS^d))$, 
then
\begin{align*}
\bE (\| \frac{1}{n} X_{n}^{T} X_{n} - \frac{1}{d+1} \cI_{d+1} \|_{\cF}^{2} )
&= 
\frac{1}{n^2} \left( \frac{d}{L +d} - \frac{d (d - 1)}{(d + L) (d + 2 L + 1) (d + 2 L -1)}\right) \binom{d + L}{d} \\
& \left ( = O(n^{- \frac{d+1}{d}})  \quad
(n \rightarrow \infty) \right ).
\end{align*}
\item[(iii)] If $\{ \x_i \}_{i = 1}^{n}$ denotes the $n$-point jittered sampling point process on $\bS^d$ 
defined in Definition~\ref{def:jspp}, then 
\[
\bE (\| \frac{1}{n} X_{n}^{T} X_{n} - \frac{1}{d+1} \cI_{d+1} \|_{\cF}^{2} ) \leq \frac{1}{n} - \frac{1}{n} \left (1 - \frac{c^2}{2 n^{2/d}} \right )^2.
\]
\end{itemize}
\end{thm}

\begin{rem}
\label{rem:1}
\begin{enumerate}
\item[(i)]
Assume that $\{ \x_i \}_{i = 1}^{n}$ denotes an $n$-point Poisson point process on $\bS^d$. 
\[
\bE (\| \frac{1}{n} X_{n}^{T} X_{n} - \frac{1}{d+1} \cI_{d+1} \|_{\cF}^{2} ) = \frac{d}{(d + 1) n}  \; \left (= O(n^{-1}) \quad  (n \rightarrow \infty)\right). 
\]
This can be easily verified based on a study by Goyal et al.~\cite{GVT98}. 
The readers interested in the generalization of this result with respect to \textit{probabilistic frames}
can refer to Ehler~\cite{E11}.
Comparing the expected value with the results of Theorem~\ref{thm:main2} reveals that 
the random point configurations generated by the DPPs converge to the FUNTF more rapidly than the Poisson point processes. 
This result occurs because the empirical covariance converges faster for the former than for the latter, 
which resonates with a lot of the literature on DPPs (e.g., Gautier et al.~\cite{GBV19}).

\item[(ii)]
The concept of $p$-frame potential is a natural generalization of frame potential (e.g., Chiristensen~\cite{C03}). 
The $p$-frame potentials of DPPs will be studied in future.
Although only three typical DPPs are addressed herein, 
various point processes on the sphere were recently proposed by 
Beltr\'an and Etayo~\cite{BE18,BE19,BE20}, 
which will be addressed in future.

\item[(iii)]
Frames have a large amount of applications in signal processing, sampling theory, wavelet theory and so on (e.g., Casazza et al.~\cite{CKP13}).
We hope that DPPs will be applied as approximate frames useful for such applications.

\end{enumerate}
\end{rem}


\section{DPPs on the Sphere}
\label{sect:dpp}

In accordance with a study by Hough et al.~\cite{HKPV09}, 
let $S$ be a locally compact Hausdorff space with a countable basis 
and $\mu$ be the Radon measure on $S$. 
Let $\cX$ be a random point process on $S$.
Then, $\cX$ is \textit{simple} if there are no coincidence points almost surely.
Letting $\cK: S^2 \rightarrow \bC$ 
be a measurable function, DPP can be defined as follows.
\begin{teigi}
\label{def:dpp}
The point process $\cX$ on $S$ is a DPP
with kernel $\cK$ if it is simple and its 
$k$-point correlation functions $\rho_k: S^k \rightarrow \R_{\geq 0}$
with respect to the measure $\mu$ satisfy
\[
\rho_k (\x_1, \ldots, \x_k) = \det (\cK (\x_i, \x_j))_{1 \leq i, j \leq k}
\]
for every $k \geq 1$; i.e.,
for any Borel function $h: S^k \rightarrow [0, \infty)$, 
we have
\begin{equation*}
\bE \big [ \sum_{\x_1, \ldots, \x_k \in \cX}^{\neq} h(\x_1, \ldots, \x_k) \big ] 
= \int_S \cdots \int_S 
\rho_k (\x_1, \ldots, \x_k) h(\x_1, \ldots, \x_k) \; d \mu(\x_1) \cdots d \mu(\x_k).
\end{equation*}
Here $\sum_{\x_1, \ldots, \x_k \in \cX}^{\neq}$ represents a multi-sum over the $n$-tuples 
of $\cX$ 
whose components are all pairwise distinct.
\end{teigi}

First, we describe the previously studied \textit{spherical ensemble}~\cite{AZ14,HKPV09,K06}, 
which can be realized, for example, using random matrices.
Let $A_n$ and $B_n$ be independent $n \times n$ random matrices with 
independent and identically distributed standard complex Gaussian entries.
The set of eigenvalues 
$\{ z_1, z_2, \ldots, z_n \}$ of $A_n^{-1} B_n$  form a DPP 
on the complex plane with kernel $\cK (z, w) = (1 + z \bar{w})^{n-1}$
with respect to the measure $\frac{n}{\pi (1 + |z|^2)^{n + 1}} \; dm(z)$, 
where $m$ denotes the Lebesgue measure on the complex plane $\bC$.
The number of eigenvalues of $A_{n}^{-1} B_{n}$ is almost surely equal to $n$.

Let $\gamma$ be the stereographic projection of the unit sphere $\bS^2$ 
from the north pole onto the plane $\{ (t_1, t_2, 0) \mid t_1, t_2 \in \R \}$.
Then, $\{\x_i = \gamma^{-1}(z_i) \mid 1 \leq i \leq n \}$ forms a DPP on $\bS^2$;
see also Figure~1.
We call such a DPP a \textit{spherical ensemble on $\bS^2$}.
For this case, Alishahi and Zamani~\cite{AZ14} showed that 
a two-point correlation function on $\bS^2$ can be given as 
\begin{equation*}
\rho_2 (\x, \y) = n^2  \left \{ 1 - \left ( 1- \frac{\|\x - \y\|^2}{4}\right )^{n-1} \right \}. 
\end{equation*}

\begin{figure}[thb]
\begin{center}
\includegraphics[width=73mm]{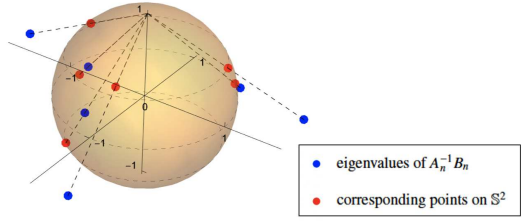} 
\end{center}
\caption{Eigenvalues of $A_n^{-1} B_n$ and their corresponding points on $\bS^2$}
\end{figure}

Next, we consider harmonic ensembles on $\bS^d$. 
Several studies have investigated such a DPP; for example, 
its Riesz energies were studied by Beltr\'an et al.~\cite{BMO16}.
Let 
\begin{equation}
\label{eq:QR}
Q_{\ell}(x) := Q_{\ell}^{d}(x) = \frac{d + 2 \ell - 1}{d - 1} C_{\ell}^{((d - 1)/2)} (x), \quad R_{L}(x) 
:= \sum_{\ell = 0}^{L} Q_{\ell}(x), 
\end{equation}
where $C_{\ell}^{(d - 1)/2)}$ is the usual Gegenbauer polynomial, i.e, 
an orthogonal polynomial on the interval  $[-1,1]$ with respect to the weight function $(1 -x^2)^{d/2 -1}$.
Note that $R_L(\langle \x, \y \rangle)$ is the $L$-th reproducing kernel for polynomial space $\cP_L(\bS^d)$, 
which is the vector space of polynomials of degree at most $L$ in $d + 1$ variables restricted to $\bS^d$. 

A DPP on $\bS^d$ associated with 
$\cK(\x, \y) = R_L(\langle \x, \y \rangle)$ is called a harmonic ensemble on $\bS^d$.  
The number of points in the harmonic ensemble is almost surely $n = \dim (\cP_L (\bS^d))
= R_L(1) = \binom{d + L}{d} + \binom{d + L -1}{d}$.
For this case,  a two-point correlation function is given as
\begin{equation*}
\rho_2 (\x, \y)  
= \det 
\begin{bmatrix}
R_{L} (1) & R_{L} (\langle \x, \y \rangle) \\
R_{L} (\langle \x, \y \rangle) & R_{L}(1) 
\end{bmatrix} 
= R_{L}(1)^2 - R_{L} (\langle \x, \y \rangle)^2.
\end{equation*}

Finally, we consider the jittered sampling point processes. 
Recently, Brauchart et al.~\cite{BSSW14} showed that the jittered sampling point process is determinantal. 
The following is a definition of the jittered sampling point process discussed here.
Let $\sigma_d$ be the normalized surface measure on $\bS^d$, 
and, for a given subset $A$ of $\bS^d$,
$\diam \; A := \sup\{ \| \x - \y \| \mid \x, \y \in A \}$ be the diameter of $A$.
\begin{teigi}
\label{def:jspp}
Let $\{ D_{i,n} \}_{i = 1}^{n}$
be an equal-area partition of $\bS^d$ into $n$ pairwise distinct subsets
with small diameters; $\bigcup_{j = 1}^{n} D_{j, n} = \bS^d$, 
where $\sigma_d (D_{j,n} \cap D_{k,n}) = 0$ for all $j, k = 1, \ldots, n$ 
with $j \neq k$ and $\sigma_{d}(D_{j,n}) = 1/n$.
Furthermore, $\diam D_{j,n} 
\leq c/n^{1/d} < \sqrt{2}$ for some $c$, independent of $n$. 
Then, letting  $\x_i := \x_{i,n}$ be a point selected randomly from $D_{i,n}$ with respect to the uniform measure on $D_{i,n}$, we call $\{ \x_{i} \}_{i = 1}^{n}$ a jittered sampling point process on $\bS^d$. 
\end{teigi}
Then, the corresponding kernel of the jittered sampling point process is 
$\cK (\x, \y) = \sum_{i = 1}^{n} \frac{\1_{D_{i,n}} (\x) \1_{D_{i,n}} (\y) }{\sigma_d (D_{i, n})}$.

\begin{rem}
In Definition~\ref{def:jspp}, 
the condition on the upper bound $\sqrt{2}$ of the diameter is a technical condition to prove Theorem~\ref{thm:main1}~(iii);
however, the existence of such partitions 
$\{ D_{i,n} \}_{i = 1}^{n}$
is guaranteed by Leopardi~\cite{L09}.
\end{rem}

\section{Proofs of Theorems~\ref{thm:main1} and \ref{thm:main2}}
\label{sect:proofs}

\subsection{Proof of Theorem~\ref{thm:main1}} 
\noindent \textbf{Proof of Theorem~\ref{thm:main1}~(i).} 
Assume that $\{ \x_i \}_{i = 1}^{n}$ represents an $n$-point spherical ensemble on $\bS^2$.
It can be easily verified that 
\begin{equation}
\label{eq:inner-distance}
|\langle \x, \y \rangle|^{2} = \left (1 - \frac{\| \x - \y \|^2}{2} \right)^{2} 
= 1 - \| \x - \y \|^2 + \frac{1}{4} \| \x - \y \|^4, \quad \x, \y \in \bS^2. 
\end{equation}

Alishahi and Zamani~\cite[Theorem 1.3(ii)]{AZ14} showed that 
for any $s \geq 0$, 
\begin{equation}
\label{eq:AZ}
 \bE ( \sum_{i \neq j} \|\x_i - \x_j\|^s ) = \frac{2^{1 + s}}{2 + s} n^2 -
 \frac{2^{s} \Gamma(n) \Gamma (1 + s/2)}{\Gamma (n + 1 + s/2)} n^2.
\end{equation}
Thus, by combining 
(\ref{eq:inner-distance}) and (\ref{eq:AZ}), 
we obtain 
\begin{align*}
\bE (\textrm{FP}(\{ \x_i \}_{i = 1}^{n}))
&= \bE ( \sum_{i \neq j}
|\langle \x_i, \x_j \rangle|^2)
= n^2 - \bE(\sum_{i \neq j}
\|\x_i - \x_j\|^2 )
+ \frac{1}{4} \bE(\sum_{i \neq j}
\|\x_i - \x_j\|^4) \\
&= n^2 - \left ( 2n^2 - \frac{4n}{n + 1} \right) + \frac{1}{4} \left ( \frac{16}{3} n^2 - \frac{32 n}{(n + 1)(n + 2)}  \right) 
=  \frac{n^2}{3} + \frac{4 n^2}{(n + 1)(n+2)}.
\end{align*}

\noindent \textbf{Proof of Theorem~\ref{thm:main1}~(ii).} 
We first recall the polynomial $Q_l$ defined in (\ref{eq:QR}).
Let $a_{l+1} = \frac{l+1}{d+2l+1}$, 
$b_{l -1} = (1 - a_{l-1})$, $b_{-1} = 0, b_0 = 1$ and $Q_{-1} = 0$.
The polynomial $Q_L$ satisfies 
\begin{align}
t Q_l (t) &= a_{l+1} Q_{l + 1} (t) + b_{l-1} Q_{l-1}(t), 
\label{eq:three-term} \\
\int_{-1}^{1} Q_l (t) Q_m (t) &(1 - t^2)^{d/2 - 1} \; dt = \delta_{l,m} \frac{|\bS^d|}{|\bS^{d-1}|} Z(d, \ell), 
\label{eq:orthogonality}
\end{align}
where $Z(d, \ell) = \dim (\textrm{Harm}_{\ell} (\R^{d + 1})) = \binom{d + \ell}{d} - \binom{d + \ell - 2}{d}$ 
and $\textrm{Harm}_{\ell} (\R^{d + 1})$ is the vector space of homogeneous and  harmonic polynomials of exact degree $\ell$ in $d + 1$ variables.

Using the three-term relation in~(\ref{eq:three-term}) and orthogonality in~(\ref{eq:orthogonality}),
the following lemma can be easily obtained.
\begin{lem}
\label{lem:1}
\begin{itemize}
\item[(i)] $\ds
\int_{-1}^{1} Q_i (t) Q_j (t) t (1 - t^2)^{d/2-1} \; dt 
=  (a_{j + 1} \delta_{i, j+1} + b_{j - 1} \delta_{i, j - 1} ) \frac{|\bS^{d}|}{|\bS^{d-1}|} Z(d, i)$.
\item[(ii)] 
$\ds \int_{-1}^{1} Q_i (t) Q_j (t) t^2 (1 - t^2)^{d/2-1} \; dt 
= (a_{i + 1} a_{j + 1} \delta_{i+1, j+1}+ a_{i + 1}  b_{j - 1} \delta_{i + 1, j-1}) \frac{|\bS^d|}{|\bS^{d-1}|} Z(d, i+1) $ \\
$\ds + (b_{i - 1}  a_{j + 1} \delta_{i-1, j+1} +  b_{i - 1} b_{j -1} \delta_{i - 1, j-1}) \frac{|\bS^d|}{|\bS^{d-1}|} Z(d, i-1).$
\end{itemize}
Here, $|\bS^d| = \frac{2 \pi^{\frac{d+1}{2}}}{\Gamma(\frac{d + 1}{2})}$ denotes the surface area of $\bS^d$.
\end{lem}

By combining the definition of $R_L$ in (\ref{eq:QR}) and Lemma~\ref{lem:1}, 
we obtain the following lemma.
\begin{lem}
\label{lem:2}
\begin{itemize}
\item[(i)] 
$\ds \int_{-1}^{1} R_L (t)^2 (1 - t^2)^{d/2 - 1} \; dt
= \frac{|\bS^{d}|}{|\bS^{d-1}|} \Big ( \binom{d + L}{d} + \binom{d + L - 1}{d} \Big )$.

\item[(ii)] 
$\ds \int_{-1}^{1} R_L (t)^2 t^2 (1 - t^2)^{d/2 - 1} \; dt
=  \frac{|\bS^d|}{|\bS^{d-1}|} \Big \{
\frac{d (L + 1)}{L (d + 2 L + 1)} + \frac{d + 4 L - 2}{d + 2 L -1} \Big \} \binom{d+L-1}{d}$.
\end{itemize}
\end{lem}

Now, we prove Theorem~\ref{thm:main1}~(ii).
Assume that $\{ \x_i \}_{i = 1}^{n}$ represents an $n$-point harmonic ensemble on $\bS^d$ with $n = \dim (\cP_L(\bS^d))$.
The continuous $s$-energy for the normalized Lebesgue measure is 
\[
V_{s}(\bS^d) = \int_{\bS^d} \int_{\bS^d} \| \x - \y \|^{-s} \; d \sigma_{d}(\x) d \sigma_{d}(\y)
= 2^{d -s -1} \frac{\Gamma(\frac{d+1}{2})\Gamma(\frac{d - s}{2})}{\sqrt{\pi}\Gamma(d - \frac{s}{2})} 
= 2^{d -s -1} \frac{\Gamma(\frac{d}{2})\Gamma(\frac{d - s}{2})}{\Gamma(d - \frac{s}{2})} \frac{|\bS^{d-1}|}{|\bS^{d}|}.
\]
Because for any positive integer $\ell$ (using the Funk-Hecke formula (e.g., M\"uller~\cite[Theorem 6]{M66})), 
\begin{align*}
\bE (\sum_{i \neq j}\|\x_i - \x_j\|^{2\ell})
&= n^2 V_{-2\ell}(\bS^d) - \int_{\bS^d} \int_{\bS^d}  R_L(\langle \x, \y \rangle)^2 \| \x - \y \|^{2\ell} \;
d \sigma_{d}(\x) d \sigma_{d}(\y) \\
&= n^2 V_{-2\ell}(\bS^d) - \frac{2^{\ell}|\bS^{d-1}|}{|\bS^d|} \int_{-1}^{1} R_L (t)^2 (1 - t)^{\ell} (1 - t^2)^{d/2 -1} \; dt, 
\end{align*}
we obtain 
\begin{align*}
\bE( \textrm{FP}(\{ \x_i \}_{i = 1}^{n})) 
& = n^2 - \bE(\sum_{i \neq j}\|\x_i - \x_j\|^2) 
+ \frac{1}{4} \bE(\sum_{i \neq j}\|\x_i - \x_j\|^4) \\
&= n^2 (1 - V_{-2}(\bS^{d}) + \frac{1}{4} V_{-4}(\bS^d)) +
\frac{|\bS^{d-1}|}{|\bS^{d}|} \int_{-1}^{1} R_L(t)^2 (1 - t^2)^{d/2 - 1} \; dt  \\ 
& \qquad - \frac{|\bS^{d-1}|}{|\bS^{d}|} \int_{-1}^{1} R_L(t)^2 t^2 (1 - t^2)^{d/2 - 1} \; dt \\
&= \frac{1}{d+1}n^2 + \left( \frac{d}{L +d} - \frac{d (d - 1)}{(d + L) (d + 2 L + 1) (d + 2 L -1)}\right) \binom{d + L}{d},
\end{align*}
where we use Lemma~\ref{lem:2} in the final equation.

\bigskip

\noindent \textbf{Proof of Theorem~\ref{thm:main1}~(iii).}
We use a previously proposed method by Brauchart et al.~\cite{BSSW14}.
Let $\{ D_{i,n} \}_{i = 1}^{n}$
be an equal-area partition of $\bS^d$ defined in Definition~\ref{def:jspp}.
Because each $D_{j,n}$ is equipped with the measure $\mu_{j, n} := \frac{\sigma_d |_{D_{j,n}}}{\sigma_d (D_{j,n})}$, 
\begin{align}
\bE (\textrm{FP} (\{ \x_i \}_{i = 1}^{n}) ) 
\nonumber 
&= \int_{D_{1, n}} \cdots \int_{D_{n,n}}  \textrm{FP} (\{ \x_i \}_{i = 1}^{n}) \; d \mu_{1, n} (\x_1) \cdots d \mu_{n, n} (\x_n) 
\nonumber \\
&= n + n^2 \int_{\bS^d} \int_{\bS^d} |\langle \x, \y \rangle|^2 d \sigma_d (\x) d 
 \sigma_d (\y) 
-  \sum_{j = 1}^{n} \int_{D_{j, n}} \int_{D_{j, n}} |\langle \x, \y \rangle|^2 
d \mu_{j, n} (\x) d \mu_{j, n} (\y)
\nonumber \\
&= \frac{n^2}{d + 1} + n 
 -  \sum_{j = 1}^{n} \int_{D_{j, n}} \int_{D_{j, n}} |\langle \x, \y \rangle|^2 
d \mu_{j, n} (\x) d \mu_{j, n} (\y).
\label{eq:expect-FP}
\end{align}

Because $\langle \x, \y \rangle^2 = (1 - \frac{1}{2} \| \x - \y \|^2)^2$ for $\x, \y \in \bS^d$
and $\diam D_{j,n} \leq c/n^{1/d} < \sqrt{2}$, 
\begin{align}
\sum_{j = 1}^{n} \int_{D_{j, n}} \int_{D_{j, n}} |\langle \x, \y \rangle|^2
d \mu_{j, n} (\x) d \mu_{j, n} (\y) 
&= \sum_{j = 1}^{n} \int_{D_{j, n}} \int_{D_{j, n}} (1 - \frac{1}{2} \| \x - \y \|^2)^{2}
d \mu_{j, n} (\x) d \mu_{j, n} (\y) \nonumber \\
&\geq \sum_{j = 1}^{n} (1 - \frac{1}{2} (\diam D_{j,n})^2)^{2}
\geq \sum_{j = 1}^{n} (1 - \frac{c^2}{2 n^{2/d}})^{2} = n (1 - \frac{c^2}{2 n^{2/d}})^{2}.
\label{eq:3rd-term}
\end{align}
Further, 
\begin{align}
\sum_{j = 1}^{n} \int_{D_{j, n}} \int_{D_{j, n}} |\langle \x, \y \rangle|^2 
d \mu_{j, n} (\x) d \mu_{j, n} (\y) 
\leq \sum_{j = 1}^{n} \int_{D_{j, n}} \int_{D_{j, n}} 
d \mu_{j, n} (\x) d \mu_{j, n} (\y) = n.
\label{eq:3rd-term2}
\end{align}

Finally, we obtain the desired result by combining 
(\ref{eq:expect-FP}), (\ref{eq:3rd-term}) and (\ref{eq:3rd-term2}).

\begin{rem}
To more precisely evaluate the frame potential,
we must obtain some tiling of the unit sphere that would enable closed-form computations instead of the bound derived in~(\ref{eq:3rd-term}), which is a future task.
\end{rem}

\subsection{Proof of Theorem~\ref{thm:main2}}

The $(i,j)$-th entry of the random matrix $X_{n}^{T} X_{n}$
associated with an $n$-point DPP $\{ \x_k \}_{k = 1}^{n}$ on $\bS^d$
is given as 
$(X_{n}^{T} X_{n})_{i, j} = \sum_{k = 1}^{n} \x_k^{(i)} \x_k^{(j)}$, 
where $\x_k = (\x_k^{(1)}, \x_{k}^{(2)}, \ldots, \x_k^{(d + 1)})^T \in \bS^d$.
For any $k, l$, 
$1 = (\sum_{i = 1}^{d+1} (\x_{k}^{(i)})^2)  (\sum_{i = 1}^{d+1} (\x_{l}^{(i)})^2)
= \sum_{1\leq i,j \leq d+1} (\x_k^{(i)} \x_l^{(j)})^2$, 
therefore, we obtain
\begin{align*}
\| \frac{1}{n} X_{n}^{T} X_{n} - \frac{1}{d + 1} \cI_{d+1} \|_{\cF}^{2} 
= \frac{1}{n^2} \textrm{FP}(\{ \x_k\}_{k = 1}^{n}) - \frac{1}{d + 1}.
\end{align*}

\noindent \textbf{Proof of Theorem~\ref{thm:main2}~(i).}
Let $\{ \x_k \}_{k = 1}^{n}$ represent an $n$-point spherical ensemble on $\bS^2$.
Using Theorem~\ref{thm:main1}~(i), we obtain
\[
\bE (\| \frac{1}{n} X_{n}^{T} X_{n} - \frac{1}{3} \cI_{3} \|_{\cF}^{2})
= \frac{1}{n^2} \bE (\textrm{FP}(\{ \x_k\}_{k = 1}^{n})) - \frac{1}{3} 
= \frac{4}{(n+1)(n+2)}.
\]

\noindent \textbf{Proof of Theorem~\ref{thm:main2}~(ii).}
Let $\{ \x_k \}_{k=1}^{n}$ represent an $n$-point harmonic ensemble on $\bS^d$ with $n = \dim(\cP_{L}(\bS^d))$.
Using Theorem~\ref{thm:main1}~(ii), we obtain 
\begin{align*}
\bE (\| \frac{1}{n} X_{n}^{T} X_{n} - \frac{1}{d+1} \cI_{d+1} \|_{\cF}^{2})
&= \frac{1}{n^2} \bE(\textrm{FP}(\{ \x_k\}_{k = 1}^{n})) - \frac{1}{d+1} \\
&= \frac{1}{n^2} \left( \frac{d}{L +d} - \frac{d (d - 1)}{(d + L) (d + 2 L + 1) (d + 2 L -1)}\right) \binom{d + L}{d}.
\end{align*}

\noindent \textbf{Proof of Theorem~\ref{thm:main2}~(iii).} 
Let $\{ \x_i \}_{i = 1}^{n}$ represent the $n$-point jittered sampling point process on $\bS^d$.
Using Theorem~\ref{thm:main1}~(iii), we obtain
\begin{align*}
\bE (\| \frac{1}{n} X_{n}^{T} X_{n} - \frac{1}{d+1} \cI_{d+1} \|_{\cF}^{2})
= \frac{1}{n^2} \bE(\textrm{FP}(\{ \x_k\}_{k = 1}^{n})) - \frac{1}{d+1} 
\leq \frac{1}{n} - \frac{1}{n} (1 - \frac{c^2}{2 n^{2/d}})^2.
\end{align*}





\noindent \textbf{Acknowledgements}

This study is partially supported by Grant-in-Aid for Young Scientists (B) 16K17645 and 
Grant-in-Aid for Scientific Research (C) 20K03736 by the Japan Society for the Promotion of Science.
The author would like to thank the anonymous referees for their valuable comments and helpful suggestions.

\end{document}